\documentclass[11pt,reqno]{amsart}

\usepackage[T1]{fontenc}
\usepackage{lmodern}
\usepackage[letterpaper,margin=1in]{geometry}
\usepackage{microtype}
\usepackage{amsmath,amssymb,amsthm,mathtools}
\usepackage{tikz-cd}
\usepackage[colorlinks=true,linkcolor=blue,citecolor=blue,urlcolor=blue]{hyperref}
\newtheorem{theorem}{Theorem}[section]
\newtheorem*{theorem*}{Theorem}
\newtheorem*{conjecture}{Conjecture}
\newtheorem{proposition}[theorem]{Proposition}
\newtheorem{lemma}[theorem]{Lemma}
\newtheorem{corollary}[theorem]{Corollary}

\newtheorem*{question*}{Question}
\theoremstyle{definition}
\newtheorem{definition}[theorem]{Definition}
\theoremstyle{remark}
\newtheorem{remark}[theorem]{Remark}

\newcommand{\ZZ}{\mathbb{Z}}
\newcommand{\QQ}{\mathbb{Q}}
\newcommand{\RR}{\mathbb{R}}
\newcommand{\CC}{\mathbb{C}}
\newcommand{\PP}{\mathbb{P}}

\newcommand{\cO}{\mathcal{O}}
\newcommand{\cA}{\mathcal A}
\newcommand{\cE}{\mathcal E}
\newcommand{\cM}{\mathcal M}
\newcommand{\cU}{\mathcal U}

\newcommand{\cZ}{\mathcal Z}

\newcommand{\HH}{\mathrm{HH}}
\newcommand{\End}{\operatorname{End}}
\newcommand{\Ext}{\operatorname{Ext}}
\newcommand{\Hom}{\operatorname{Hom}}
\newcommand{\id}{\operatorname{id}}
\newcommand{\tr}{\operatorname{tr}}
\newcommand{\pr}{\operatorname{pr}}
\newcommand{\rk}{\operatorname{rk}}
\newcommand{\ch}{\operatorname{ch}}
\newcommand{\td}{\operatorname{td}}
\newcommand{\Ku}{\operatorname{Ku}}
\newcommand{\Db}{\mathrm D^{\mathrm b}}

\newcommand{\RHom}{\mathrm R\mathcal Hom}
\newcommand{\rHom}{\mathrm R\!\operatorname{Hom}}

\newcommand{\sEnd}{\mathcal{E}\!nd}
\newcommand{\Pic}{\operatorname{Pic}}

\newcommand{\Stab}{\operatorname{Stab}}
\newcommand{\tors}{\mathrm{tors}}
\newcommand{\rO}{\mathrm O}

\numberwithin{equation}{section}

\title[Derived Torelli for Kuznetsov components]
{Derived Torelli for Kuznetsov components of cubic fourfolds}
\author{Hanfei Guo}
\address{Tongji University, School of Mathematical Sciences, 1239 Siping Road, Shanghai 200092, China}
\email{hanfeiguo23@gmail.com, hfguo@tongji.edu.cn}

\author{Ruxuan Zhang}
\address{Fudan University, School of Mathematical Sciences,
220 Handan Road, Shanghai 200433, China}
\email{rxzhang18@fudan.edu.cn}

\date{}

\begin{document}

\begin{abstract}
We prove that every orientation-preserving Hodge isometry between
the integral Mukai lattices of the Kuznetsov components of two smooth
cubic fourfolds is induced by an exact equivalence.
The proof uses Markman--Mehrotra's comonad reconstruction method
to recover both Kuznetsov components from a common Fourier--Mukai
comonad on a moduli space of stable objects.
As an application, we prove that the Kuznetsov components of two smooth
cubic fourfolds are equivalent whenever their Fano varieties of lines
are derived equivalent.
\end{abstract}

\maketitle
\tableofcontents

\section{Introduction}\label{sec:introduction}

Throughout, we work over the complex numbers $\CC$.
All equivalences are assumed to be $\CC$-linear.

\subsection{Derived Torelli theorem}
Let $S$ be a projective K3 surface. Its Mukai lattice
\[
  \widetilde H(S,\ZZ)
  :=H^0(S,\ZZ)\oplus H^2(S,\ZZ)\oplus H^4(S,\ZZ)
\]
carries the Mukai pairing and a weight-two Hodge structure. By the
derived Torelli theorem of Mukai and Orlov, this Hodge lattice determines
the bounded derived category of $S$.

\begin{theorem*}[Mukai--Orlov {\cite{Mukai87,Orlov97}}]
Let $S$ and $S'$ be projective K3 surfaces. Then $\Db(S)\simeq\Db(S')$ if
and only if there is a Hodge isometry
$\widetilde H(S,\ZZ)\xrightarrow{\ \sim\ }\widetilde H(S',\ZZ)$.
\end{theorem*}

Cubic fourfolds provide another source of K3 categories. For a smooth
cubic fourfold $X\subset\PP^5$, consider the semiorthogonal component
\cite{Kuznetsov10}
\[
  \Ku(X):=\langle\cO_X,\cO_X(1),\cO_X(2)\rangle^\perp\subset\Db(X),
\]
which is a noncommutative K3 surface: its Serre functor is the shift
$[2]$, and its Hochschild homology agrees with that of a K3 surface.
By \cite[\S2]{AddingtonThomas14}, $\Ku(X)$ has an integral Mukai lattice
$\widetilde H(\Ku(X),\ZZ)$ equipped with a weight-two Hodge structure
and abstractly isometric to the Mukai lattice of a K3 surface. This
lattice carries a natural orientation, which we recall in
Section~\ref{sec:mukai-lattice}.

Let $X$ and $Y$ be smooth cubic fourfolds.
By \cite[Proposition~3.4]{Huybrechts17} and
\cite[Theorem~1.3]{LiPertusiZhaoHearts23},
any exact equivalence $\Ku(X)\simeq\Ku(Y)$ induces a Hodge isometry
$\widetilde H(\Ku(X),\ZZ)\simeq\widetilde H(\Ku(Y),\ZZ)$.
This leads to the following question of Huybrechts; see also
\cite[Question~5.4]{LiPertusiZhao18},
\cite[Question~3.24]{MacriStellari19},
\cite[Conjecture~6.1]{HuybrechtsUpdate25}, and
\cite[Conjecture~2.12]{BohningBothmerMarquand26}.

\begin{question*}[Huybrechts]
Suppose that there is an
orientation-preserving Hodge isometry
\[
  \widetilde H(\Ku(X),\ZZ)\xrightarrow{\ \sim\ }\widetilde H(\Ku(Y),\ZZ).
\]
Are $\Ku(X)$
and $\Ku(Y)$ derived equivalent?
\end{question*}

By \cite[Theorem~1.5]{Huybrechts17}, the answer is affirmative for
nonspecial cubic fourfolds and for very general members of each
Hassett divisor. Using the stability conditions constructed in
\cite{BLMSStability23}, one also obtains an affirmative answer
when the algebraic part of
$\widetilde H(\Ku(X),\ZZ)$ contains a nonzero isotropic class: in this
case $\Ku(X)$ is equivalent to the
derived category of a twisted K3 surface
\cite[Proposition~33.1]{BLMNPSFamilies19}, and the claim follows from
the twisted derived Torelli theorem of Huybrechts and Stellari
\cite[Theorem~0.1]{HuybrechtsStellari07}.

We prove the following theorem, which answers
Huybrechts' question in general.

\begin{theorem}\label{thm:torelli}
Let $X$ and $Y$ be smooth cubic fourfolds. If there is an
orientation-preserving Hodge isometry
\[
  g\colon\widetilde H(\Ku(X),\ZZ)\xrightarrow{\ \sim\ }
  \widetilde H(\Ku(Y),\ZZ),
\]
then there is an exact equivalence
$\Psi\colon\Ku(X)\xrightarrow{\sim}\Ku(Y)$ such that $\Psi_*=g$ on the
integral Mukai lattices.
\end{theorem}

By \cite[Lemma~2.3]{Huybrechts17}, the orientation condition can be omitted
if one asks only for an equivalence $\Psi$, without requiring $\Psi_*=g$
on the integral Mukai lattice.
This proves \cite[Conjecture~2.13]{BohningBothmerMarquand26}.

\subsection{Fano varieties of lines}
Let $X$ and $Y$ be smooth cubic fourfolds, and let $F(X)$ and $F(Y)$
be their Fano varieties of lines. Both are hyper-K\"ahler fourfolds of
$\mathrm{K3}^{[2]}$-type \cite{BeauvilleDonagi85}. By
\cite[Theorem~A]{KemboiSegal25}, $\Db(F(X))$ is equivalent to the
categorical symmetric square of $\Ku(X)$. In particular, an equivalence
$\Ku(X)\simeq\Ku(Y)$ induces an equivalence
$\Db(F(X))\simeq\Db(F(Y))$. As an application of
Theorem~\ref{thm:torelli}, we prove the converse
in general, which answers a question of Brooke, Frei and Marquand
\cite[Introduction]{BrookeFreiMarquand24}; see also
\cite[Introduction]{MeinsmaMoschetti25}.

\begin{corollary}\label{cor:fano-ku}
Let $X$ and $Y$ be smooth cubic fourfolds. If
$\Db(F(X))\simeq\Db(F(Y))$, then there is an exact equivalence
$\Ku(X)\simeq\Ku(Y)$.
\end{corollary}

When $\Ku(X)$ is equivalent to the derived
category of an untwisted K3 surface, Corollary~\ref{cor:fano-ku} follows from
\cite[Theorem~1]{AddingtonRationality16} and
\cite[Corollaries~9.6--9.7]{Beckmann23}.

\subsection{Strategy of the proof}\label{sec:strategy}
In Section~\ref{sec:moduli}, we use the Hodge isometry $g$ to construct
a common fine moduli space $M\cong F(Y)$ for the two Kuznetsov components.
We establish the local structure, stability and base-change properties
of the relative Ext sheaves in Section~\ref{sec:ext-sheaf}.
In Section~\ref{sec:deformation}, we identify the two relative Ext sheaves
after suitable choices of universal families, using deformation to
moduli spaces of sheaves on K3 surfaces and the rigidity of the modular sheaf in
\cite{MMVRigidity17}.
We then compare the associated Fourier--Mukai comonads and split
the units in Section~\ref{sec:torelli-proof} to reconstruct
an equivalence inducing $g$ on the Mukai lattices.
Finally, in Section~\ref{sec:fano-application}, we prove
Corollary~\ref{cor:fano-ku} by lattice cancellation.

\medskip
\noindent\textbf{Acknowledgements.}
We are grateful to Chunyi Li for many comments and suggestions
on the first draft of this paper. We thank Lie Fu, Lisa Marquand and
Alex Perry for useful discussions on this project. We also thank Simone
Bille, Zhiyuan Li, Zhiyu Liu, Laura Pertusi and Qizheng Yin for helpful conversations.

\medskip
\noindent\textbf{AI disclosure.}
We had long considered proving Theorem~\ref{thm:torelli} using
Markman--Mehrotra's comonad reconstruction method
\cite[Proposition~1.4]{MarkmanMehrotraTransforms15}. The second author
proposed applying this strategy to the problem in a talk at Shanghai Jiao
Tong University in December 2025. We began the project then, but had not yet
overcome the technical difficulties.

In August 2026, we shared this strategy with GPT~6 Astra. After extensive
interaction, it produced a complete proof of Theorem~\ref{thm:torelli}.
Our work with LLMs consisted mainly of reading, understanding, checking
and revising the proof.

After we circulated the draft among some members of the community,
an anonymous expert conducted an independent test, prompting Astra to
explore comonad reconstruction and the construction of fine moduli
spaces. He concluded that these heuristic suggestions alone
would not readily enable Astra to reproduce the result, at least
within one day. We therefore believe that the work is
worth posting on arXiv.

We wrote the paper and take responsibility for its correctness.

\section{Moduli spaces of stable objects}\label{sec:moduli}

In this section, we first recall the Mukai lattices and stability
conditions of the Kuznetsov components of cubic fourfolds. We then use a
Hodge isometry between the Mukai lattices of two such components to
construct a common moduli space. A similar strategy was discussed in
\cite[\S5.3]{LiPertusiZhao18}.

\subsection{Mukai lattices}\label{sec:mukai-lattice}
Let $X\subset\mathbb{P}^5$ be a smooth cubic fourfold. The derived category
of $X$ admits the semiorthogonal decomposition
\cite{Kuznetsov10}
\[
  \Db(X)=\langle\Ku(X),\cO_X,\cO_X(1),\cO_X(2)\rangle.
\]
Let $\pr_X\colon\Db(X)\to\Ku(X)$ be the projection functor,
that is, the left adjoint of the inclusion
$\Ku(X)\hookrightarrow\Db(X)$. The category $\Ku(X)$ is
two-Calabi--Yau: its Serre functor is $[2]$.

The Mukai lattice is
\[
  \widetilde H(\Ku(X),\ZZ)
  :=\{a\in K_{\mathrm{top}}(X):\chi([\cO_X(i)],a)=0\text{ for }i=0,1,2\},
\]
with pairing $(a,b):=-\chi(a,b)$, where $K_{\mathrm{top}}(X)$ denotes
complex topological $K$-theory. It is even and unimodular of signature
$(4,20)$ and carries a weight-two Hodge structure
\cite[\S2.3]{AddingtonThomas14}.

The \emph{algebraic Mukai lattice} of $\Ku(X)$
\[
  N(\Ku(X)):=\widetilde H(\Ku(X),\ZZ)\cap\widetilde H^{1,1}(\Ku(X))
\]
is primitive and has signature $(2,\rk N(\Ku(X))-2)$.
The class map identifies the numerical Grothendieck group of $\Ku(X)$
with $N(\Ku(X))$ by \cite[Proposition~2.4]{AddingtonThomas14}.
For $E\in\Ku(X)$, let $[E]\in N(\Ku(X))$ denote its Mukai class.
For a line $l\subset X$, define
$\lambda_{i,X}:=[\pr_X\cO_l(i)]$ for $i=1,2$. These classes generate the
primitive sublattice
$A_2(X):=\ZZ\lambda_{1,X}\oplus\ZZ\lambda_{2,X}$, whose intersection pairing is
\[
  \begin{pmatrix}2&-1\\-1&2\end{pmatrix};
\]
see \cite[\S2.4]{AddingtonThomas14}, where the opposite sign convention
is used. By \cite[Proposition~2.3]{AddingtonThomas14}, the Mukai vector
induces a Hodge isometry
\[
  A_2(X)^\perp\xrightarrow{\ \sim\ }H^4_{\mathrm{prim}}(X,\ZZ)(1),
\]
where the right-hand side carries the negative intersection pairing.

For a generator $\zeta_X$ of $\widetilde H^{2,0}(\Ku(X))$, the classes
$\operatorname{Re}\zeta_X,\operatorname{Im}\zeta_X,\lambda_{1,X},\lambda_{2,X}$
form an oriented basis of a positive definite four-plane in
$\widetilde H(\Ku(X),\RR)$. The orientation is independent of
$\zeta_X$ and induces an orientation on every positive definite
four-plane by orthogonal projection.

\begin{definition}
Let $X, Y$ be smooth cubic fourfolds. A Hodge isometry
\[
  \widetilde H(\Ku(X),\ZZ)\xrightarrow{\ \sim\ }\widetilde H(\Ku(Y),\ZZ)
\]
is \emph{orientation-preserving} if it preserves these orientations.
\end{definition}

Such an isometry need not map the sublattice $A_2(X)$ to $A_2(Y)$.

\subsection{Stability conditions}\label{sec:stability}
Let $\Stab^\dagger(\Ku(X))$ be the connected component of the space of
numerical Bridgeland stability conditions on $\Ku(X)$ containing those
constructed in \cite{BLMSStability23}. A stability condition
$\sigma\in\Stab^\dagger(\Ku(X))$ has central charge
$Z_\sigma(-)=(\Omega_\sigma,-)$ for a unique
$\Omega_\sigma\in N(\Ku(X))\otimes\CC$. Let
\[
  \mathcal P_0^+(X)\subset N(\Ku(X))\otimes\CC
\]
be the set of classes $\Omega$ such that
$(\operatorname{Re}\Omega,\operatorname{Im}\Omega)$ is an oriented basis
of a positive definite two-plane, with the orientation determined by
that of Section~\ref{sec:mukai-lattice} and the Hodge plane, and such
that $(\Omega,\delta)\neq0$ for every $\delta\in N(\Ku(X))$ with
$\delta^2=-2$.
By \cite[Theorem~29.1]{BLMNPSFamilies19}, the map
$\sigma\mapsto\Omega_\sigma$ restricts to a covering over
$\mathcal P_0^+(X)$.

For a fixed class $v\in N(\Ku(X))$, the walls are locally finite,
and the stable objects of class $v$ are unchanged within each chamber
\cite[Lemma~12.13]{BLMNPSFamilies19}.
We say that $\sigma$ is \emph{$v$-generic} if it does not lie on a wall
for $v$; when $v$ is primitive, stability and semistability then coincide.
For primitive $v\in N(\Ku(X))$ with $v^2=2$
and $v$-generic $\sigma\in\Stab^\dagger(\Ku(X))$, the coarse moduli
space $M_\sigma(v)$ of stable objects of class $v$ and fixed phase
is a projective hyper-K\"ahler fourfold of
$\mathrm{K3}^{[2]}$-type, and the Mukai morphism
\begin{equation}\label{eq:mukai-morphism}
  \theta_v\colon v^\perp\xrightarrow{\ \sim\ }H^2(M_\sigma(v),\ZZ)
\end{equation}
is a Hodge isometry for the Beauville--Bogomolov--Fujiki form
\cite[Theorem~29.2]{BLMNPSFamilies19}; here $v^\perp$ is taken in
the full Mukai lattice $\widetilde H(\Ku(X),\ZZ)$.
We use the dual determinant convention of
\cite[Definition~5.4 and Remark~5.5]{BayerMacriProjectivity12} for $\theta_v$.
The Bayer--Macr\`i class of $\sigma$ on $M_\sigma(v)$ is
\begin{equation}\label{eq:bm-class}
  \ell_\sigma:=\theta_v\bigl(\operatorname{Im}(-\Omega_\sigma/Z_\sigma(v))\bigr)
  \in N^1(M_\sigma(v))_\RR.
\end{equation}
The class is nef by \cite[Theorem~21.25]{BLMNPSFamilies19}.

At a point $[E]\in M_\sigma(v)$, the Kodaira--Spencer map identifies
$T_{[E]}M_\sigma(v)$ with $\Ext^1(E,E)$; the symplectic form is the
Yoneda pairing followed by the Serre trace
\cite[\S31]{BLMNPSFamilies19}.

\subsection{The common moduli space}\label{sec:common-moduli}
From now on, we fix an orientation-preserving Hodge isometry
\[
  g\colon\widetilde H(\Ku(X),\ZZ)\xrightarrow{\ \sim\ }
  \widetilde H(\Ku(Y),\ZZ).
\]
Set
\[
  \begin{aligned}
    v_Y&:=\lambda_{1,Y}+\lambda_{2,Y},&\qquad u_Y&:=\lambda_{1,Y},\\
    v_X&:=g^{-1}(v_Y),& u_X&:=g^{-1}(u_Y).
  \end{aligned}
\]
Then $v_X^2=u_X^2=2$ and $(u_X,v_X)=1$; in particular, $v_X$ is primitive. By
\cite{LiPertusiZhao18} and \cite[Example~29.6]{BLMNPSFamilies19}, there
is a $v_Y$-generic stability condition $\sigma_Y\in\Stab^\dagger(\Ku(Y))$
with $M_{\sigma_Y}(v_Y)\cong F(Y)$. We identify their second cohomology
groups via this isomorphism.

\begin{proposition}\label{prop:common-moduli}
The stability condition $\sigma_Y$ can be chosen so that there exist
a $v_X$-generic stability condition
$\sigma_X\in\Stab^\dagger(\Ku(X))$ and an isomorphism
\begin{equation}\label{eq:common-moduli}
  f\colon M_{\sigma_X}(v_X)\xrightarrow{\ \sim\ }M_{\sigma_Y}(v_Y)\cong F(Y)
\end{equation}
such that $(f^{-1})^*=h$, where
\begin{equation}\label{eq:moduli-isometry}
  h:=\theta_{v_Y}\circ g\circ\theta_{v_X}^{-1}\colon
  H^2(M_{\sigma_X}(v_X),\ZZ)\xrightarrow{\ \sim\ }H^2(F(Y),\ZZ).
\end{equation}
Moreover, both moduli spaces are fine.
\end{proposition}

\begin{proof}
Let $Z_Y$ be the central charge of $\sigma_Y$. Moving $\sigma_Y$ within
its chamber, we may assume that
\begin{equation}\label{eq:generic-charge}
  \operatorname{Im}\frac{Z_Y(a)}{Z_Y(v_Y)}\neq0
  \qquad\text{for all }a\in N(\Ku(Y))\setminus\QQ v_Y,
\end{equation}
since the excluded central charges form a countable union of proper
real hypersurfaces. As $g$ identifies $\mathcal P_0^+(X)$ with
$\mathcal P_0^+(Y)$, we may choose
$\sigma_X\in\Stab^\dagger(\Ku(X))$ with central charge
$Z_X:=Z_Y\circ g$ by \cite[Theorem~29.1]{BLMNPSFamilies19}.
By \eqref{eq:generic-charge}, $\sigma_X$ is $v_X$-generic.

The Hodge isometry $h$ satisfies
$h(\ell_{\sigma_X})=\ell_{\sigma_Y}$ by \eqref{eq:bm-class}.
Deforming $\sigma_X$ within its chamber, we obtain a neighborhood
of $\ell_{\sigma_X}$ in $N^1(M_{\sigma_X}(v_X))_\RR$ consisting of nef classes
\cite[Theorem~12.11, Lemma~12.13, and Theorem~21.25]{BLMNPSFamilies19}.
Thus $\ell_{\sigma_X}$ is ample, and the same argument applies to
$\ell_{\sigma_Y}$. By global Torelli
and the monodromy description for manifolds of $\mathrm{K3}^{[2]}$-type
\cite[Theorem~1.3 and Lemma~9.2]{MarkmanSurvey11}, we obtain an
isomorphism $f$ as in \eqref{eq:common-moduli} with $(f^{-1})^*=h$.

Finally, since $\chi(u_X,v_X)=\chi(u_Y,v_Y)=-1$, both moduli spaces
are fine; see \cite[Theorem~4.6.5]{HuybrechtsLehn10} and
\cite[proof of Lemma~33.2]{BLMNPSFamilies19}.
\end{proof}

Let $M$ be the common moduli space in \eqref{eq:common-moduli}, with
universal families $\cU_X$ on $X\times M$ and $\cU_Y$ on $Y\times M$.
Their Mukai morphisms satisfy
\begin{equation}\label{eq:marking-compatibility}
  \theta_{v_Y}\circ g|_{v_X^\perp}=\theta_{v_X}.
\end{equation}

\section{The relative Ext sheaf}\label{sec:ext-sheaf}
We establish the local structure and base-change properties of
relative Ext sheaves by local calculations in the analytic topology, generalizing
\cite[Proposition~4.1]{MarkmanCharacteristic11}.
We then deduce their slope stability
from \cite[Theorem~1.4]{MarkmanNaturality16}.

By faithful flatness of analytic local rings over algebraic local rings
\cite[\S2, no.~6, Corollary~1]{SerreGAGA56}, local freeness, reflexivity
and flatness of algebraic coherent sheaves can be checked after
analytification, as can the invertibility of algebraic morphisms between
them. We use the same notation for analytified objects in local calculations.

\subsection{Local structure along the diagonal}\label{sec:local-model}
Let $Z$ be a smooth projective variety and let $\cA\subset\Db(Z)$ be
an admissible triangulated subcategory whose Serre functor $S_\cA$
is equipped with an isomorphism $S_\cA\cong[2]$.
Let $M$ be a smooth connected algebraic moduli space
of dimension four parametrizing pairwise nonisomorphic stable objects
of a fixed numerical class and phase for a Bridgeland stability
condition on $\cA$. Let $\cU$ be a universal
complex on $Z\times M$ with bounded coherent cohomology and fibres in
$\cA$, possibly twisted by the pullback of an $\cO^\times$-valued
cocycle $\alpha$ on $M$.
We assume that the Kodaira--Spencer map
\begin{equation}\label{eq:ks}
  T_mM\xrightarrow{\sim}\Ext^1_\cA(\cU_m,\cU_m)
\end{equation}
is an isomorphism for every $m\in M$.

Let $p_i\colon M\times M\to M$ be the projections, let $\Delta$ be
the diagonal, and set $\cU_i:=(\id_Z\times p_i)^*\cU$.
For the projection $\pi\colon Z\times M\times M\to M\times M$, let
$\RHom_\pi(-,-):=\mathrm R\pi_*\RHom(-,-)$, and set
\[
  \mathcal R:=\RHom_\pi(\cU_1,\cU_2),\qquad E:=\mathcal H^1(\mathcal R).
\]
The perfect complex $\mathcal R$ is twisted by
$p_1^*\alpha^{-1}p_2^*\alpha$, and its derived fibre at $(m,n)$ is
$\rHom_\cA(\cU_m,\cU_n)$ by relative Hom base change
\cite[Lemma~3.21]{BLMNPSFamilies19}.
By stability, Serre duality and \eqref{eq:ks}, the self-Ext dimensions
are $1,4,1$ in degrees $0,1,2$, so the common numerical class has
Euler pairing $-2$. Hence
\begin{equation}\label{eq:ext-dimensions}
  \dim\Ext^i_\cA(\cU_m,\cU_n)=
  \begin{cases}
    1,\,4,\,1 & \text{for }i=0,1,2\text{ if }m=n,\\
    0,\,2,\,0 & \text{for }i=0,1,2\text{ if }m\ne n,
  \end{cases}
\end{equation}
with all other Ext groups zero. The Yoneda pairing and the Serre trace
induce a nondegenerate holomorphic two-form $\omega_M$ through
\eqref{eq:ks}. When $M$ is hyper-K\"ahler, this form is symplectic.

Let $\beta\colon\operatorname{Bl}_\Delta(M\times M)\to M\times M$
be the blow-up, with exceptional divisor $D=\PP(T_M)$, projection
$p\colon D\to M$, and tautological line
$\ell:=\cO_D(-1)\subset p^*T_M$.
Orthogonal complements in $p^*T_M$ are taken with respect to
$p^*\omega_M$.

\begin{proposition}\label{prop:local-model}
With the notation above, the following hold.
\begin{enumerate}
\item[(a)] The trace identifies $\mathcal H^2(\mathcal R)$ with $\cO_\Delta$,
and $\mathcal H^0(\mathcal R)=0$. The sheaf $E$ is reflexive of rank two and
locally free away from $\Delta$.
\item[(b)] The twisted sheaf
\begin{equation}\label{eq:normalized-pullback}
  V:=\bigl(\beta^*E/\tors\bigr)(D)
\end{equation}
is locally free of rank two.
\item[(c)] The twisting of $V$ is canonically trivial along $D$, and
\begin{equation}\label{eq:exceptional-restriction}
  V|_D\cong\ell^\perp/\ell.
\end{equation}
\item[(d)] If $M$ is a compact irreducible hyper-K\"ahler fourfold, then $E$ is
slope-stable for $p_1^*\omega+p_2^*\omega$ for every K\"ahler class
$\omega$ on $M$, and $\End(E)=\CC$.
\end{enumerate}
\end{proposition}

\begin{proof}
Away from $\Delta$, we have $\mathcal R\cong E[-1]$ with $E$ locally free
of rank two. We work in a small analytic polydisc about $(m,m)$ and
trivialize the twisting. Let $R$ be the regular analytic local ring at
$(m,m)$ and $I$ the diagonal ideal.
By \eqref{eq:ext-dimensions}, we may choose a minimal free representative
\begin{equation}\label{eq:minimal-complex}
  \mathcal R_{\mathrm{loc}}^\bullet:=[R\xrightarrow{x}R^4\xrightarrow{b}R]
\end{equation}
in degrees $0,1,2$. The differentials vanish modulo $I$, since the
derived fibres along the diagonal have constant dimensions.
The identity, the Kodaira--Spencer map and the trace identify its
three terms modulo $I$ with $\cO_M,T_M,\cO_M$.

For a tangent vector $\xi=(\xi_1,\xi_2)$, the first-order term of
the Hom differential on a homogeneous element $f$ is
$\xi_2f-(-1)^{|f|}f\xi_1$. Hence
\[
  \partial_\xi x(\id)=\xi_2-\xi_1,\qquad
  \tr\bigl(\partial_\xi b(f)\bigr)=\tr\bigl((\xi_2-\xi_1)f\bigr),
\]
where the second identity uses graded cyclicity of the Serre trace
\cite[p.~94]{CaldararuWillerton10}.
The entries of $x$ therefore generate $I$ and form a regular sequence.
By the Koszul relations for $bx=0$, we obtain
\begin{equation}\label{eq:alternating-matrix}
  b=x^{\mathsf t}A,\qquad A^{\mathsf t}=-A.
\end{equation}
By the second identity in the first-order calculation, $A$ modulo the
maximal ideal represents the symplectic pairing. Thus $A$ is invertible, and the
entries of $b$ also generate $I$.

Consequently, $\mathcal H^0(\mathcal R)=0$ and $\mathcal H^2(\mathcal R)=R/I$;
the local identifications of $\mathcal H^2(\mathcal R)$ glue via the trace.
Set $K_0:=\ker b$. We prove that $E$ is reflexive using the sequences
\[
  0\longrightarrow K_0\longrightarrow R^4\xrightarrow{b}I
    \longrightarrow0,\qquad
  0\longrightarrow R\xrightarrow{x}K_0\longrightarrow E
    \longrightarrow0
\]
and the depth lemma. At a prime $\mathfrak p\supset I$
of height $d\ge4$, we obtain
\[
  \operatorname{depth}(K_0)_{\mathfrak p}\ge d-2,\qquad
  \operatorname{depth}E_{\mathfrak p}\ge d-2\ge2.
\]
Combining these inequalities with local freeness away from $\Delta$,
we conclude that $E$ is torsion-free and satisfies Serre's condition $S_2$,
hence is reflexive \cite[Tag~0AVT]{StacksProject}.

On a blow-up chart, let $x=tu$, where $t$ is an equation of $D$ and
$u$ is a unimodular column. Set $q:=u^{\mathsf t}A$. Then $b=tq$,
$qu=0$, and $q$ is unimodular, so
\[
  Q:=\ker(q\colon\cO^4\to\cO)/\cO u
\]
is locally free of rank two. Let $e_1,\ldots,e_4$ be the standard
basis of $R^4$. The Koszul generators $b_i e_j-b_j e_i$
of $K_0$ pull back to $t(q_i e_j-q_j e_i)$, which span $t\ker q$.
Thus there is a surjection $\beta^*E\to tQ$ that is an isomorphism
away from $D$. Its kernel is exactly the torsion subsheaf, so
\[
  \beta^*E/\tors\cong Q(-D),\qquad V\cong Q.
\]
Along $D$, the vector $u$ spans the tautological line in $p^*T_M$,
and $q$ is contraction with this line under the symplectic pairing.
This proves \eqref{eq:exceptional-restriction} locally.
The maps $u$ and $q$ along $D$ are the first normal differentials of
the algebraic complex $\mathrm L\beta^*\mathcal R$, under the
Kodaira--Spencer and trace identifications.
On overlaps, comparison maps between the local complexes induce the
Kodaira--Spencer identifications in degree one modulo $I$. The induced
maps are unchanged under homotopy, since the differentials vanish modulo
$I$. Thus the local isomorphisms \eqref{eq:exceptional-restriction}
glue to an algebraic isomorphism. The twisting is canonically trivial along $D$
because $p_1=p_2$ on $\Delta$.

If $M$ is compact irreducible hyper-K\"ahler, we verify
\cite[Assumption~1.3]{MarkmanNaturality16} using (a)--(c).
The slope stability in (d) follows from
\cite[Theorem~1.4]{MarkmanNaturality16}. Since a slope-stable twisted
sheaf is simple \cite[proof of Lemma~4.1.9]{LieblichOrbifold11},
we have $\End(E)=\CC$.
\end{proof}

\subsection{Families}\label{sec:families}
Let $C$ be a smooth connected algebraic curve. A \emph{family of moduli
fourfolds over $C$} consists of a smooth projective morphism
$\mathcal V\to C$, admissible triangulated subcategories
$\cA_c\subset\Db(\mathcal V_c)$ for $c\in C$, a smooth proper morphism
$q\colon\cM\to C$ with connected four-dimensional fibres, and a
universal complex $\cU$ on $\mathcal V\times_C\cM$ with bounded
coherent cohomology. This complex may be twisted by an
$\cO^\times$-valued cocycle pulled back from $\cM$; its fibres over
$\cM_c$ lie in $\cA_c$.
Let $p_i\colon\cM\times_C\cM\to\cM$ be the projections, and let
$\cU_i$ be the corresponding pullbacks of $\cU$ to
$\mathcal V\times_C\cM\times_C\cM$. Let $\pi$ be the projection to
$\cM\times_C\cM$.
Assume that every fibre satisfies the hypotheses of
Section~\ref{sec:local-model} and that the Serre traces are compatible
in families. By \cite[Lemma~3.21]{BLMNPSFamilies19}, the perfect complex
\[
  \mathcal R:=\RHom_\pi(\cU_1,\cU_2)
  \quad\text{on}\quad\cZ:=\cM\times_C\cM
\]
commutes with arbitrary base change. Let $\Delta\subset\cZ$ be the
relative diagonal and set $\cE:=\mathcal H^1(\mathcal R)$.

To compare the relative Ext sheaves along a deformation in
Section~\ref{sec:ext-comparison}, we need them to form a flat family
with the expected fibres.

\begin{lemma}\label{lem:base-change}
For a family of moduli fourfolds as above, the sheaf $\cE$ is flat
over $C$. For every $c\in C$, the natural map
\[
  \cE|_{\cM_c\times\cM_c}
  \longrightarrow
  \mathcal H^1\RHom_{\pi_c}(\cU_{c,1},\cU_{c,2})
\]
is an isomorphism, where $\pi_c$ and $\cU_{c,i}$ are the restrictions
over $c$.
\end{lemma}

\begin{proof}
We apply the local calculation of Proposition~\ref{prop:local-model}
over $C$, following
\cite[Proposition~4.1]{MarkmanCharacteristic11}.
In the analytic local ring $R$ at a point of $\Delta$, after trivializing
the twisting, we have
\[
  \mathcal R\cong[R\xrightarrow{x}R^4\xrightarrow{x^{\mathsf t}A}R]
\]
in degrees $0,1,2$, with $A$ invertible and alternating. By the
fibrewise Kodaira--Spencer isomorphism and Nakayama's lemma, the
entries of $x$ generate the relative diagonal ideal. The same holds
for the entries of $x^{\mathsf t}A$. Thus the trace identifies
$\mathcal H^2(\mathcal R)$ with $\cO_\Delta$, which is flat over $C$.

Fix $c\in C$ and choose a local parameter $t$ at $c$.
Let $i_c\colon\cZ_c\hookrightarrow\cZ$ be the inclusion. By
derived base change, $\mathrm Li_c^*\mathcal R\cong
\RHom_{\pi_c}(\cU_{c,1},\cU_{c,2})$, whose degree-zero
cohomology vanishes by Proposition~\ref{prop:local-model}(a).
The cohomology sequence of
\[
  \mathcal R\xrightarrow{t}\mathcal R
  \longrightarrow i_{c*}\mathrm Li_c^*\mathcal R
\]
therefore shows that multiplication by $t$ on $\cE$ is injective.
Since $C$ is a smooth curve, $\cE$ is flat over $C$.
Multiplication by $t$ on $\mathcal H^2(\mathcal R)$ is also injective,
so the same sequence identifies
$\cE/t\cE$ with $i_{c*}\mathcal H^1(\mathrm Li_c^*\mathcal R)$.
This proves the fibre identification.
\end{proof}

Replacing $\cU$ by $\cU\otimes L$, for a line bundle $L$ on $\cM$,
replaces $\cE$ by $\cE\otimes p_1^*L^{-1}\otimes p_2^*L$.
This follows from the projection formula for relative derived Hom.

\section{Comparison with deformed modular sheaves}\label{sec:deformation}

In this section, we compare the relative Ext sheaves associated with the universal
families $\cU_X$ and $\cU_Y$ over the common moduli space $M$ of
Section~\ref{sec:common-moduli}. Using the intrinsic endomorphism algebra
of deformed modular sheaves \cite{MMVRigidity17}, we choose the universal families
so that their relative Ext sheaves are isomorphic.

\subsection{Stable and rigid families}
\begin{lemma}\label{lem:relative-isomorphism-locus}
Let $q\colon Z\to C$ be a smooth proper holomorphic map to a connected
smooth complex curve, with connected compact K\"ahler fibres. Let $V_1,V_2$
be $C$-flat twisted coherent sheaves with the same
$\cO^\times$-valued cocycle.
Suppose their fibres are reflexive, have the same rank and determinant,
and are slope-stable for a common K\"ahler class on each fibre. If
\[
  V_{1,c_*}\cong V_{2,c_*},\qquad
  \Ext^1(V_{1,c_*},V_{1,c_*})=0
\]
for some $c_*\in C$, then $V_{1,c}\cong V_{2,c}$ for every $c\in C$.
\end{lemma}

\begin{proof}
On a small analytic disc in $C$, the perfect complex
$\mathrm Rq_*\RHom(V_1,V_2)$
admits a finite locally free representative whose derived fibres compute
the Ext groups of the fibres, by analytic cohomology and base change
\cite[\S\S6--7]{GrauertModuli60}.
A nonzero map between the two stable fibres has full rank. Since
their determinants agree, the map is invertible in codimension one;
by reflexivity, its inverse is defined on the whole fibre
\cite[Lemma~2.9]{MarkmanMehrotraTorelli13}. Hence the
isomorphism locus is
\[
  I:=\{c\in C:\dim\Hom(V_{1,c},V_{2,c})\ge1\},
\]
a closed analytic subset by semicontinuity. Since
$\Ext^1$ vanishes at $c_*$, we can lift the fibre isomorphism to a
morphism on a small analytic disc about $c_*$ by cohomology and base change.
This morphism is an isomorphism after shrinking the disc. Thus $I$ contains
an open set, and $I=C$ by the identity principle.
\end{proof}

\subsection{The canonical comparison}\label{sec:ext-comparison}
We first recall the modular sheaf on a moduli space of sheaves on a
K3 surface. Let $S$ be a projective K3 surface, $v$ a primitive Mukai
vector with $v^2=2$, and $H$ a $v$-generic polarization such that
$N:=M_H(v)$ is nonempty. Choose a possibly twisted universal sheaf
$\cU$ on $S\times N$. Let
$\pi_N\colon S\times N\times N\to N\times N$ and
$q_i\colon S\times N\times N\to S\times N$ be the projections.
The \emph{geometric modular sheaf} is the relative Ext sheaf
\[
  E_N:=\mathcal H^1\RHom_{\pi_N}(q_1^*\cU,q_2^*\cU);
\]
see \cite[\S1.2, equation~(1.3)]{MMVRigidity17}.
It is reflexive of rank two and locally free away from the diagonal,
with fibre $\Ext^1_S(F,G)$ at a pair $([F],[G])$ with $F\not\cong G$
\cite[Proposition~4.1]{MarkmanCharacteristic11}.
If $\cU$ is twisted by a Brauer class $\alpha$ pulled back from $N$,
then $E_N$ is $p_1^*\alpha^{-1}p_2^*\alpha$-twisted, where
$p_i\colon N\times N\to N$ are the projections.

Now let $M$ be an arbitrary manifold of $\mathrm{K3}^{[2]}$-type.
Choose markings of $N$ and $M$ in the same connected component of the
moduli space of marked manifolds, and a twistor path $\gamma$ joining
them. Deforming $E_N$ along the corresponding diagonal twistor
families, we obtain a possibly twisted rank-two reflexive sheaf
$E_\gamma$ on $M\times M$, called a \emph{deformed modular sheaf}
\cite[\S1.2 and Theorem~1.9]{MMVRigidity17}.
The sheaf is considered up to tensoring by line bundles. By
\cite[Theorem~1.11(2)--(4)]{MMVRigidity17}, the unordered pair
\[
  \{\sEnd(E_\gamma),\sEnd(E_\gamma^\vee)\}
\]
of isomorphism classes of algebras is independent of the initial
moduli space, the markings and the twistor path.
For a rank-two reflexive sheaf,
$E_\gamma^\vee\cong E_\gamma\otimes\det(E_\gamma)^{-1}$, so these
two algebras are isomorphic. Denote their common isomorphism class
by $\cA_M$. Thus $\cA_M$ is intrinsic to $M$, although we construct
it using a twistor deformation of a geometric modular sheaf.

For manifolds of nonmaximal Picard number,
\cite[Theorem~1.11]{MarkmanMehrotraTorelli13} holds unconditionally
by \cite[\S1.2 and Theorem~1.11(1)]{MMVRigidity17}.
Let $\cM\to B$ be an analytic deformation disc of $M$ with
$\rho(M)<21$. After shrinking $B$, there is a $B$-flat family
$\mathcal F$ of rank-two reflexive twisted sheaves on
$\cM\times_B\cM$ with $\sEnd(\mathcal F_0)\cong\cA_M$
\cite[Proposition~2.14, Lemmas~4.5 and~4.8, and Corollary~4.13]{MarkmanMehrotraTorelli13}.
For $\rho(\cM_b)<21$, the fibre $\mathcal F_b$ is a deformed
modular sheaf up to a line bundle twist by
\cite[Lemma~2.15, Lemma~4.11 and Theorem~1.11(3)(b)]{MarkmanMehrotraTorelli13} and
\cite[\S1.2 and Theorem~1.11]{MMVRigidity17}.
Thus
\[
  \sEnd(\mathcal F_b)\cong\cA_{\cM_b},\qquad
  \Ext^1_{\cM_b\times\cM_b}(\mathcal F_b,\mathcal F_b)=0.
\]
We use the vanishing \cite[Theorem~1.11(1)]{MMVRigidity17} in
Lemma~\ref{lem:relative-isomorphism-locus} to compare $\mathcal F$ with $\cE$.

\begin{proposition}\label{prop:canonical-comparison}
Let $\cM\to C$ be a family of moduli fourfolds of
$\mathrm{K3}^{[2]}$-type as in Section~\ref{sec:families}, and let
$0,1\in C$. Assume that $\cM_1$ is a moduli space of Gieseker-stable
sheaves on a projective K3 surface with geometric modular sheaf $\cE_1$.
If there is a path $\gamma$ from $1$ to $0$ in $C$ with
$\rho(\cM_c)<21$ for all $c\in\gamma$, then the relative Ext sheaf $\cE$ satisfies
\[
  \sEnd(\cE_0)\cong\cA_{\cM_0}.
\]
\end{proposition}

\begin{proof}
We compare $\cE$ with the local families $\mathcal F$ constructed above by checking
the hypotheses of Lemma~\ref{lem:relative-isomorphism-locus} on
successive small analytic discs along $\gamma$.

On such a disc $B$, set $Z:=\cM\times_B\cM$. On a sufficiently fine
common analytic cover, trivialize the local determinants of $\cE|_Z$
and $\mathcal F$ and choose square roots of the determinants of their
gluing maps. For each sheaf, replace its gluing maps by
\[
  \varphi_{ij}':=(\det\varphi_{ij})^{-1/2}\varphi_{ij}.
\]
Denote the resulting twisted sheaves by $\cE'$ and $\mathcal F'$.
They have trivial determinant, $\mu_2$-valued twisting cocycles and
the same endomorphism algebras as $\cE|_Z$ and $\mathcal F$, respectively, by
\cite[Proposition~2.14]{MarkmanMehrotraTorelli13} and
\cite[\S2.1]{MarkmanCharacteristic11}.

Their twisting classes in $H^2(Z,\mu_2)$ are determined by one fibre,
since the relative diagonal has codimension four
\cite[\S1.3]{MMVRigidity17} and $Z\to B$ is differentiably trivial
\cite[\S1]{EhresmannConnections52}.
Thus, if $\sEnd(\cE_{c_*})\cong\cA_{\cM_{c_*}}$ at some $c_*\in B$,
these classes coincide. Multiplying the gluing maps of $\mathcal F'$
by suitable signs, we can choose the same twisting cocycle for both
sheaves. On that fibre, $\cE'_{c_*}$ and $\mathcal F'_{c_*}$ differ
by a line bundle $L$ with
$L^{\otimes2}\cong\cO$
\cite[Proposition~2.14]{MarkmanMehrotraTorelli13}.
Since $\Pic(\cM_{c_*}\times\cM_{c_*})$ has no torsion,
$\cE'_{c_*}\cong\mathcal F'_{c_*}$.

Both $\cE'$ and $\mathcal F'$ are flat with rank-two reflexive fibres,
by Lemma~\ref{lem:base-change} and the construction of $\mathcal F$.
Their fibres are stable for a common diagonal K\"ahler class by
\cite[Theorem~1.4]{MarkmanNaturality16} and
Proposition~\ref{prop:local-model}(d), and
$\Ext^1(\mathcal F'_{c_*},\mathcal F'_{c_*})=0$ by the preceding
discussion. Hence Lemma~\ref{lem:relative-isomorphism-locus} applies
on $B$. Starting with the geometric modular sheaf $\cE_1$, we apply
the lemma on overlapping discs along $\gamma$ to obtain
$\sEnd(\cE_0)\cong\cA_{\cM_0}$.
\end{proof}

\subsection{Deformation to moduli spaces on K3 surfaces}\label{sec:geometric-deformation}
Recall the common moduli space
\[
  M=M_{\sigma_X}(v_X)=M_{\sigma_Y}(v_Y)
\]
of Section~\ref{sec:common-moduli}. For $W\in\{X,Y\}$, set
$N_W:=N(\Ku(W))$. By \eqref{eq:mukai-morphism},
\[
  \rho(M)=\rk N_W-1.
\]
If $\rho(M)=21$, then $N_W$ has signature $(2,20)$.
Its orthogonal complement in the unimodular Mukai lattice
$\widetilde H(\Ku(W),\ZZ)$ has rank two,
so the discriminant group of $N_W$ has at most two generators
\cite[Proposition~1.6.1]{Nikulin80}. Hence $N_W$ contains a
hyperbolic plane $U$ by \cite[Corollary~1.13.5]{Nikulin80}.
By \cite[Lemma~32.5 and Proposition~33.1]{BLMNPSFamilies19}, there are
equivalences $\Phi_W\colon\Ku(W)\xrightarrow{\sim}\Db(S_W)$ for
projective K3 surfaces $S_W$, compatible with distinguished stability
components. The Hodge isometry
\[
  (\Phi_Y)_*\circ g\circ(\Phi_X)_*^{-1}\colon
  \widetilde H(S_X,\ZZ)\xrightarrow{\ \sim\ }\widetilde H(S_Y,\ZZ)
\]
is orientation-preserving. By
\cite[Theorem~0.1]{HuybrechtsStellari07}, it is induced by a
Fourier--Mukai equivalence $\Phi$.
Thus $\Phi_Y^{-1}\circ\Phi\circ\Phi_X$ induces $g$, proving
Theorem~\ref{thm:torelli} in this case.

For the rest of this subsection, assume that $\rho(M)<21$.
For each moduli description, we apply
Proposition~\ref{prop:canonical-comparison} to the relative Ext sheaf
in a family containing $M$ and a moduli space of sheaves on a K3 surface.

\begin{lemma}\label{lem:geometric-family}
For each of the two moduli descriptions of $M$ above, there is a
family $\cM\to C$ as in Section~\ref{sec:families}, with
$0,1\in C$ and $\cM_0=M$, such that $\cM_1$ is a moduli space of
Gieseker-stable sheaves on an untwisted projective K3 surface.
The fibre $\cE_1$ is its geometric modular sheaf.
There is a path $\gamma$ from $0$ to $1$ in $C$ with
$\rho(\cM_c)<21$ for all $c\in\gamma$.
\end{lemma}

\begin{proof}
We prove the assertion for $M=M_{\sigma_X}(v_X)$; an analogous argument
holds for $M=M_{\sigma_Y}(v_Y)$.
Let $N_0$ be the saturation of $\langle A_2(X),u_X,v_X\rangle$
in $\widetilde H(\Ku(X),\ZZ)$. Its rank $r$ is at most four and
its signature is $(2,r-2)$. The complement $N_0^\perp$ has signature
$(2,22-r)$ and discriminant length at most $r$, since the Mukai lattice
is unimodular \cite[Proposition~1.6.1]{Nikulin80}. Its rank is at least
twenty, so it contains a
hyperbolic plane $U=\ZZ e\oplus\ZZ f$
\cite[Corollary~1.13.5]{Nikulin80}. Set
\[
  t:=e-7f,\qquad w:=3u_X-v_X+t,\qquad N_1:=N_0\oplus\ZZ t.
\]
Then $t^2=-14$, $w^2=0$ and $(v_X,w)=1$. The lattice $N_1$ is
primitive: pairing with $e,f$ forces the coefficient of $t$ in
any integral saturation class to be integral.
Moreover,
\[
  N_1\cap A_2(X)^\perp
  =(N_0\cap A_2(X)^\perp)\oplus\ZZ t
\]
contains no new classes of squares $-2,-6$: a vector
with nonzero $t$-coefficient has square at most $-14$.
Choose a very general period orthogonal to $N_1$ in the component
of the $N_0$-Hodge locus containing $X$. By the smoothness of $X$ and
the preceding calculation, this period avoids the hyperplanes excluded
by the cubic period theorem \cite[Theorem~1.1]{LazaPeriod10}.
Hence it corresponds to a smooth cubic $X'$ with
$N(\Ku(X'))=N_1$. Since $w$ is isotropic of divisibility one, we apply
\cite[Lemma~32.5 and Proposition~33.1]{BLMNPSFamilies19} to identify
$\Ku(X')$ with $\Db(S)$ for an untwisted projective K3 surface $S$,
compatibly with distinguished stability components.

By \cite[Corollary~32.1 and Lemma~9.8]{BLMNPSFamilies19}, there is
a family $\cM\to C$ as in Section~\ref{sec:families}, with
$\cM_0=M$ and $\cM_1$ a moduli space of stable objects in $\Ku(X')$.
By \cite[Proposition~3.5]{BudurZhangFormality19}, a Fourier--Mukai
equivalence identifies $\cM_1$ with a moduli space of Gieseker-stable
sheaves on $S$ and $\cE_1$ with its geometric modular sheaf.

By \eqref{eq:mukai-morphism},
$\rho(\cM_1)=\rk N_1-1=r\le4$.
If the period map is nonconstant, the points with maximal Picard
number form a countable set \cite[\S1.2]{MMVRigidity17}.
We can choose a path from $0$ to $1$ in $C$
that avoids this set. If the period map is constant, all fibres have
nonmaximal Picard number.
\end{proof}

For $W\in\{X,Y\}$, let
$\pi_W\colon W\times M\times M\to M\times M$ be the projection,
$\cU_{W,i}:=(\id_W\times p_i)^*\cU_W$, and set
\[
  E_W:=\mathcal H^1\RHom_{\pi_W}(\cU_{W,2},\cU_{W,1}).
\]
The results of Section~\ref{sec:ext-sheaf} apply to $E_W$ after
interchanging the two factors.

\begin{proposition}\label{prop:ext-comparison}
The universal families $\cU_X$ and $\cU_Y$ can be chosen so that
\begin{equation}\label{eq:ext-comparison}
  E_X\cong E_Y.
\end{equation}
\end{proposition}

\begin{proof}
Let $\tau$ be the involution of $M\times M$ interchanging the factors.
Applying Proposition~\ref{prop:canonical-comparison} to the families in
Lemma~\ref{lem:geometric-family}, separately for $X$ and $Y$, and
pulling back by $\tau$, we obtain
\[
  \sEnd(E_X)\cong\tau^*\cA_M\cong\sEnd(E_Y).
\]
Both sheaves are untwisted. By
\cite[Proposition~2.14]{MarkmanMehrotraTorelli13}, there is a line
bundle $L_{12}$ on $M\times M$ such that
$E_X\cong E_Y\otimes L_{12}$.
Since $H^1(M,\ZZ)=H^1(M,\cO_M)=0$, the K\"unneth formula and the
exponential sequence identify
$\Pic(M\times M)=p_1^*\Pic(M)\oplus p_2^*\Pic(M)$;
hence $L_{12}=L_1\boxtimes L_2$.
The two Serre forms are nonzero elements of the one-dimensional space
$H^0(M,\Omega_M^2)$, so their orthogonal complements of $\ell$ agree.
By Proposition~\ref{prop:local-model}(c), both exceptional restrictions
on the blow-up are therefore isomorphic to $Q:=\ell^\perp/\ell$, and
\[
  Q\cong Q\otimes p^*(L_1\otimes L_2).
\]
Taking determinants and using the injectivity of
$p^*\colon\Pic(M)\to\Pic(D)$, we obtain
$(L_1\otimes L_2)^{\otimes2}\cong\cO_M$.
The Picard group is torsion-free by the exponential sequence, so
$L_2=L_1^{-1}$.
Let $p_M\colon Y\times M\to M$ be the projection. The object
$\cU_Y\otimes p_M^*L_1$ is again a universal family for
$M_{\sigma_Y}(v_Y)$. Replacing $\cU_Y$ by this family replaces $E_Y$ by
$E_Y\otimes(L_1\boxtimes L_1^{-1})$ and proves
\eqref{eq:ext-comparison}.
\end{proof}

From now on, we fix universal families satisfying
\eqref{eq:ext-comparison}. This normalization preserves $\theta_{v_Y}$
on $v_Y^\perp$: for $x\in v_Y^\perp$, the determinant class defining
$\theta_{v_Y}(x)$ changes by a multiple of
$\chi(x,v_Y)=-(x,v_Y)$, which vanishes. In particular,
\eqref{eq:marking-compatibility} still holds.

\section{Proof of Theorem~\ref{thm:torelli}}\label{sec:torelli-proof}
Let $X$ and $Y$ be smooth cubic fourfolds with an
orientation-preserving Hodge isometry
\[
  g\colon\widetilde H(\Ku(X),\ZZ)\xrightarrow{\ \sim\ }
  \widetilde H(\Ku(Y),\ZZ).
\]
By the reduction in Section~\ref{sec:geometric-deformation}, we may
assume that $\rho(M)<21$.

We first use Proposition~\ref{prop:ext-comparison} to identify the
counital kernels. We then apply Addington's splitting result on $\Ku(Y)$
to compare the comonad structures. In Section~\ref{sec:reconstruction},
we split the unit on $\Ku(X)$ using Hochschild homology and
apply Markman--Mehrotra's reconstruction method to obtain a
Fourier--Mukai equivalence $\Ku(X)\simeq\Ku(Y)$. We then prove that
this equivalence induces $g$ on the Mukai lattices.
For the definitions and properties of comonads and their comodules, see
\cite[\S2]{ElaginEquivariant14} and
\cite[\S1.3]{MarkmanMehrotraTransforms15}.

\subsection{Counital kernels}\label{sec:kernel}
For $W\in\{X,Y\}$, let
\[
  G_W\colon\Db(M)\longrightarrow\Ku(W),\qquad G_W(k_m)\cong\cU_{W,m},
\]
be the Fourier--Mukai functor with kernel $\cU_W$, where $k_m$ is the
skyscraper sheaf at $m\in M$ and $\cU_{W,m}$ is the derived restriction
of $\cU_W$ to $W\times\{m\}$. Let $F_W$ be its left adjoint, and let
$i_W\colon\Ku(W)\hookrightarrow\Db(W)$ be the inclusion.
By Serre duality in $\Ku(W)$, we have
\[
  F_W\cong\Phi_{\cU_W^\vee[2]}\circ i_W.
\]
By \cite[Example~2.3]{ElaginEquivariant14} and
\cite[\S1.3]{MarkmanMehrotraTransforms15}, the adjunction
$F_W\dashv G_W$ defines a comonad $F_W G_W$ on $\Db(M)$, with
counit $\epsilon_W\colon F_W G_W\to\id_{\Db(M)}$ and comultiplication
$d_W:=F_W\eta_W G_W$, where
$\eta_W\colon\id_{\Ku(W)}\to G_W F_W$ is the unit.
The Fourier--Mukai kernel of $F_W G_W$ on $M\times M$ is
\begin{equation}\label{eq:kernel-fiber}
  K_W\cong\RHom_{\pi_W}(\cU_{W,2},\cU_{W,1})[2],
\end{equation}
whose derived fibre at $(m,n)$ is $\rHom_{\Ku(W)}(\cU_{W,n},\cU_{W,m})[2]$.
By Proposition~\ref{prop:local-model}(a), the complex $K_W$ has
exactly two cohomology sheaves, namely
\[
  \mathcal H^{-1}(K_W)=E_W,\qquad
  \mathcal H^0(K_W)\cong\cO_\Delta.
\]
We also denote by
$\epsilon_W\colon K_W\to\cO_\Delta$ the kernel morphism representing
the counit; it induces the degree-zero identification.
The truncation triangle is therefore
\begin{equation}\label{eq:kernel-triangle}
  E_W[1]\longrightarrow K_W
  \xrightarrow{\ \epsilon_W\ }\cO_\Delta
  \xrightarrow{\ \kappa_W\ }E_W[2].
\end{equation}
Here $\kappa_W\in\Ext^2_{M\times M}(\cO_\Delta,E_W)$ denotes the
connecting class.

By the local model \eqref{eq:minimal-complex}, the derived restriction
of $K_W[-2]$ to $\Delta$ has cohomology in degrees $0,1,2$.
By Grothendieck duality for the regular
embedding $\Delta\subset M\times M$ of codimension four, we have
$\Ext^j(\cO_\Delta,K_W)=0$ for $j<2$. Applying
$\Hom(\cO_\Delta,-)$ to \eqref{eq:kernel-triangle}, we find that
composition with $\kappa_W$ induces isomorphisms
\[
  \Ext^j(\cO_\Delta,\cO_\Delta)\xrightarrow{\ \sim\ }
  \Ext^{j+2}(\cO_\Delta,E_W)\qquad(j\leq0).
\]
Thus $\Ext^1(\cO_\Delta,E_W)=0$ and
$\Ext^2(\cO_\Delta,E_W)=\CC\kappa_W$ with $\kappa_W\ne0$.
After rescaling the identification in \eqref{eq:ext-comparison}, the two
connecting classes agree. Taking their shifted cones, we identify the
counital kernels:
\begin{equation}\label{eq:common-kernel}
  (K,\epsilon):=(K_X,\epsilon_X)\cong(K_Y,\epsilon_Y).
\end{equation}
Set $E:=E_X\cong E_Y$ and denote the common connecting class by $\kappa$.

\subsection{Isomorphism of comonads}\label{subsec:splitting}\label{subsec:comonad}
We prove that the counit $\epsilon$ determines the comultiplication on
$K$ uniquely. For this, we first establish a splitting of $G_YF_Y$ on
$\Ku(Y)$.

Let $W$ be a smooth cubic fourfold and let $M:=M_\sigma(v)$ be a smooth
projective fine moduli space of $\sigma$-stable objects in $\Ku(W)$ of
character $v\in N(\Ku(W))$, with $\dim M=2n\geq 2$. Let
$\Phi\colon\Ku(W)\to\Db(M)$ be the left adjoint of the Fourier--Mukai
transform $\Db(M)\to\Ku(W)$ defined by the universal family, and let
$R$ be the right adjoint of $\Phi$.
We generalize Addington's splitting conjecture
\cite{AddingtonNewDerived16} as follows.

\begin{conjecture}[Addington]
There is a decomposition
\[
  R\Phi\cong\bigoplus_{i=0}^{n-1}\id_{\Ku(W)}[-2i]
\]
under which the unit $\eta\colon\id_{\Ku(W)}\to R\Phi$ is the
inclusion of the summand with $i=0$.
\end{conjecture}

For moduli spaces of sheaves on a K3 surface, the decomposition holds in
the Hilbert scheme case, and the unit splits for general
primitive $v$ \cite[Theorem~2.2]{MarkmanMehrotraTransforms15}.
We establish the splitting for $M=M_{\sigma_Y}(v_Y)$.

\begin{proposition}\label{prop:universal-splitting}
There is an isomorphism
\begin{equation}\label{eq:addington-monad}
  G_YF_Y\cong\id_{\Ku(Y)}\oplus[-2]
\end{equation}
under which the unit $\eta_Y$ is the inclusion of the first summand.
\end{proposition}

\begin{proof}
By \eqref{eq:common-moduli}, we have $M=M_{\sigma_Y}(v_Y)\cong F(Y)$. Let
\[
  \mathcal P:=\{(y,l)\in Y\times M:y\in l\}
\]
be the universal line, with projections $q\colon\mathcal P\to Y$ and
$p\colon\mathcal P\to M$. For the inclusion
$\iota\colon\Ku(Y)\hookrightarrow\Db(Y)$, define
\[
  F_{\mathrm{Add}}:=\mathrm Rp_*\circ\mathrm Lq^*\circ\iota
  \colon\Ku(Y)\longrightarrow\Db(M),
\]
and let $R_{\mathrm{Add}}$ be its right adjoint. By
\cite[Theorem~5.1 and (5.4)]{AddingtonNewDerived16}, the triangle
\[
  \id_{\Ku(Y)}\xrightarrow{\ \eta_{\mathrm{Add}}\ }
  R_{\mathrm{Add}}F_{\mathrm{Add}}\longrightarrow[-2]
\]
splits at the level of Fourier--Mukai kernels.

Let $\mathbf L_Q$ and $\mathbf R_Q$ denote left and right mutation
through an exceptional object $Q$. For a line $l\subset Y$, let $I_l$
be its ideal sheaf. By \cite[\S5.1]{AddingtonNewDerived16},
\[
  R_{\mathrm{Add}}(k_l)\cong\mathbf L_{\cO_Y}(I_l(1)).
\]
Our universal family parametrizes the stable objects
$P_l:=\mathbf R_{\cO_Y(-1)}I_l$ of class $v_Y$, with triangle
\[
  \cO_Y(-1)[1]\longrightarrow P_l\longrightarrow I_l;
\]
see \cite[\S4]{LiPertusiZhao18}. Let
$\mathsf R:=\mathbf L_{\cO_Y}\circ(-\otimes\cO_Y(1))$ be the rotation
autoequivalence of $\Ku(Y)$
\cite[Appendix~A]{BLMSStability23}.
Tensoring the triangle by $\cO_Y(1)$ and applying $\mathbf L_{\cO_Y}$,
we obtain
\[
  \mathsf R(P_l)\cong\mathbf L_{\cO_Y}(I_l(1))
  \cong R_{\mathrm{Add}}(k_l).
\]

The two families have isomorphic simple fibres and no negative
self-extensions. By \cite[Corollary~4.3.3]{LieblichComplexes06}, the
stack of such complexes is a $\mathbb G_m$-gerbe over an algebraic
space. Since $M$ is reduced and the fibres agree, the two families
induce the same morphism to this space.
Thus, for some line bundle $L_{\mathrm{Add}}$ on $M$, we have kernel
isomorphisms
\begin{equation}\label{eq:addington-conjugation}
  \mathsf R\circ G_Y\cong R_{\mathrm{Add}}\circ T_L,
  \qquad
  F_Y\cong T_L^{-1}\circ F_{\mathrm{Add}}\circ\mathsf R,
\end{equation}
where $T_L:=(-\otimes L_{\mathrm{Add}})$.
We obtain the second isomorphism by taking left adjoints in the first.
These identifications respect the units. The line-bundle twists
cancel in the composite, and hence
\[
  G_YF_Y\cong\mathsf R^{-1}\circ
  (R_{\mathrm{Add}}F_{\mathrm{Add}})\circ\mathsf R
  \cong\id_{\Ku(Y)}\oplus[-2].
\]
The splitting of $\eta_Y$ follows from that of $\eta_{\mathrm{Add}}$.
\end{proof}

\begin{proposition}\label{prop:comonad}
The comonads $F_X G_X$ and $F_Y G_Y$ on $\Db(M)$ are isomorphic through
an isomorphism of kernels respecting the counits and comultiplications.
\end{proposition}

\begin{proof}
For kernels on $M\times M$, denote convolution by $*$, with unit
$\cO_\Delta$. Composing the splitting in
Proposition~\ref{prop:universal-splitting} on the left with $F_Y$ and
on the right with $G_Y$, we obtain
\begin{equation}\label{eq:convolution-square}
  K*K\cong K\oplus K[-2].
\end{equation}

Since $K$ has cohomology sheaves only in degrees $-1$ and $0$, we have
$\Hom(K,K[-2])=0$. By \eqref{eq:kernel-triangle}, we have
$\Hom(K,\cO_\Delta)=\CC\epsilon$ and an exact sequence
\[
  0\longrightarrow\Hom(K,E[1])\longrightarrow\End(E)
  \longrightarrow\Ext^2(\cO_\Delta,E).
\]
By Proposition~\ref{prop:local-model}(d), $\End(E)=\CC$, and the
last map sends $\id_E$ to $\kappa\ne0$. Thus $\Hom(K,E[1])=0$,
and composition with $\epsilon$ identifies $\End(K)$ with $\CC\epsilon$.
Using \eqref{eq:convolution-square}, we deduce
\[
  \Hom(K,K*K)\cong\End(K)\oplus\Hom(K,K[-2])=\CC.
\]
Under \eqref{eq:common-kernel}, let $d_W\colon K\to K*K$ represent
$F_W\eta_WG_W$ for $W\in\{X,Y\}$. Since
$(\epsilon*K)\circ d_Y=\id_K$, the map $d_Y$ spans this space.
Thus $d_X=a\,d_Y$. Applying the counit law for $d_X$, we obtain
$\id_K=a\,\id_K$, so $a=1$. Hence \eqref{eq:common-kernel}
is an isomorphism of comonads
\cite[Definition~2.2]{ElaginEquivariant14}.
\end{proof}

\subsection{Reconstruction of the equivalence}\label{sec:reconstruction}\label{subsec:torelli-conclusion}
We first split $\eta_X$, then reconstruct the equivalence and determine
its action on the Mukai lattice. In this subsection, all morphisms of
functors are induced by morphisms of Fourier--Mukai kernels. For admissible
subcategories, we use the corresponding kernels after projection.
By Proposition~\ref{prop:universal-splitting}, projection onto the
identity summand is a retraction of $\eta_Y$. To split $\eta_X$, we apply
the trace argument in
\cite[proof of Theorem~2.2(2), \S\S2.4--2.5]{MarkmanMehrotraTransforms15}
in the following form.

\begin{lemma}\label{lem:splitting-criterion}
Let $W$ be a smooth cubic fourfold, and let
$F\colon\Ku(W)\to\Db(M)$ be a Fourier--Mukai functor with right
adjoint $G$. If
\[
  F_{\HH}\colon\HH_{-2}(\Ku(W))\longrightarrow\HH_{-2}(M)
\]
is nonzero, then the unit $\eta\colon\id_{\Ku(W)}\to GF$ admits a
natural retraction (a left inverse) induced by a morphism of
Fourier--Mukai kernels.
\end{lemma}

\begin{proof}
Let $\cA:=\Ku(W)$. By additivity of Hochschild homology and the HKR
isomorphism
\cite[Theorem~4.5, Proposition~4.6 and Corollary~7.5]{KuznetsovHH09},
we have
\[
  \HH_{-2}(\cA)\cong\HH_{-2}(W)\cong H^{1,3}(W)\cong\CC.
\]
Thus $F_{\HH}$ is injective.
Let $\mathcal P_W\in\Db(W\times W)$ be the kernel of
$i_W\pr_W$, where $\pr_W$ is the left adjoint of the inclusion $i_W$.
This kernel represents the identity on $\cA$.
Fix an isomorphism $S_{\cA}\cong[2]$ at the level of these kernels.
By Serre duality, we have
\[
  \HH_{-2}(\cA)^\vee
  \cong\Hom(\id_{\cA},S_{\cA}[-2])
  \cong\End(\mathcal P_W).
\]
For $\nu\in\HH_{-2}(M)^\vee=\Hom(\id_{\Db(M)},S_M[-2])$, define
\[
  s_\nu\colon GF\xrightarrow{G\nu F}GS_M[-2]F
  \cong F^LF\xrightarrow{\epsilon^L}\id,
\]
where $F^L$ is the left adjoint of $F$ with counit $\epsilon^L$, and
$GS_M\cong S_{\cA}F^L\cong F^L[2]$.
By the trace formula
\cite[(2.16)--(2.17)]{MarkmanMehrotraTransforms15}, we obtain
\begin{equation}\label{eq:hh-trace-transpose}
  F_{\HH}^\vee(\nu)=s_\nu\circ\eta\in\End(\mathcal P_W).
\end{equation}
Since $F_{\HH}^\vee$ is surjective, we can choose $\nu$ such that
$s_\nu\eta=\id_{\mathcal P_W}$.
\end{proof}

\begin{lemma}\label{lem:split-unit}
The unit $\eta_X\colon\id_{\Ku(X)}\to G_XF_X$ admits a natural
retraction $s_X\colon G_XF_X\to\id_{\Ku(X)}$ induced by a morphism of
Fourier--Mukai kernels.
\end{lemma}

\begin{proof}
Let $A_W:=(F_W)_{\HH}$ and $B_W:=(G_W)_{\HH}$ denote the induced maps
on Hochschild homology in degree $-2$. By additivity and functoriality
of Hochschild homology and \eqref{eq:addington-monad}, we have
\[
  B_YA_Y=2\,\id_{\HH_{-2}(\Ku(Y))},
\]
since the induced map of the even shift $[-2]$ on Hochschild homology
is the identity.
In particular, $A_Y\ne0$. By \eqref{eq:common-kernel}, we also have
$A_XB_X=A_YB_Y$. Hence
\begin{equation}\label{eq:hh-nonzero}
  A_XB_XA_Y=A_YB_YA_Y=2A_Y\ne0.
\end{equation}
Thus $A_X\ne0$, and Lemma~\ref{lem:splitting-criterion} applies.
\end{proof}

We can now reconstruct the equivalence from the comonad isomorphism in
Proposition~\ref{prop:comonad}.

\begin{proposition}\label{prop:reconstruction}
There is an exact equivalence $\Psi\colon\Ku(X)\xrightarrow{\sim}\Ku(Y)$
induced by a Fourier--Mukai kernel, with kernel isomorphisms
\begin{equation}\label{eq:reconstruction-kernels}
  F_Y\Psi\cong F_X,\qquad \Psi G_X\cong G_Y.
\end{equation}
\end{proposition}

\begin{proof}
For $W\in\{X,Y\}$, set $\mathsf K_W:=F_WG_W$. By
Proposition~\ref{prop:universal-splitting} and Lemma~\ref{lem:split-unit},
both units split at the level of kernels.
Since $\Ku(W)$ is idempotent complete
\cite[Corollary~2.10]{BalmerSchlichting01}, the comparison
functor is an equivalence \cite[Proposition~2.8]{ElaginEquivariant14}:
\[
  Q_W\colon\Ku(W)\xrightarrow{\sim}
  \operatorname{Comod}(\mathsf K_W)(\Db(M)),\qquad
  x\longmapsto(F_Wx,F_W\eta_{W,x}),
\]
where $\operatorname{Comod}(\mathsf K_W)(\Db(M))$ is the category of
comodules over the comonad $\mathsf K_W$ on $\Db(M)$
\cite[Definition~2.5]{ElaginEquivariant14}.
Let $a\colon\mathsf K_X\xrightarrow{\sim}\mathsf K_Y$ denote the
comonad isomorphism in Proposition~\ref{prop:comonad}.
By \cite[Lemma~2.21]{ElaginDescent12}, it induces an equivalence $a_*$
of comodule categories via $(A,\rho)\mapsto(A,a_A\circ\rho)$.
The equivalence $a_*$ commutes with the forgetful functors.

We construct a kernel representing $Q_Y^{-1}a_*Q_X$.
Set $H:=G_YF_X$. The $\mathsf K_Y$-comodule structure on $F_X$
is defined by the morphism
\[
  \rho:=(aF_X)\circ(F_X\eta_X)\colon F_X\longrightarrow F_YH.
\]
Let $s_Y\colon G_YF_Y\to\id_{\Ku(Y)}$ be the projection onto the
identity summand in Proposition~\ref{prop:universal-splitting}, and set
\[
  e:=s_{Y,H}\circ G_Y\rho\colon H\longrightarrow H,\qquad
  p:=F_Ye,\qquad r:=\epsilon_{Y,F_X}\circ p.
\]
Here $s_{Y,H}:=s_YH$ and $\epsilon_{Y,F_X}:=\epsilon_YF_X$, with
$p\colon F_YH\to F_YH$ and $r\colon F_YH\to F_X$.
By naturality and the defining identities for a comodule, we have
\[
  e^2=e,\qquad p\rho=\rho,\qquad \rho r=p,\qquad r\rho=\id_{F_X}.
\]
These are identities of kernel morphisms.
Let $\mathcal H\in\Db(X\times Y)$ be the convolution kernel of
$i_YH\pr_X$ \cite[Lemma~57.8.3]{StacksProject}.
The idempotent $e$ on $\mathcal H$ splits by
\cite[Corollary~2.10]{BalmerSchlichting01}.
Let $\mathcal P$ be its image, with
\[
  j\colon\mathcal P\longrightarrow\mathcal H,\qquad
  q\colon\mathcal H\longrightarrow\mathcal P,\qquad
  qj=\id_{\mathcal P},\qquad jq=e.
\]
Then $\Psi:=\Phi_{\mathcal P}|_{\Ku(X)}$ defines an exact functor
from $\Ku(X)$ to $\Ku(Y)$. The maps
\[
  F_Yq\circ\rho\colon F_X\longrightarrow F_Y\Psi,\qquad
  r\circ F_Yj\colon F_Y\Psi\longrightarrow F_X
\]
are inverse kernel morphisms respecting the comodule structures. Hence
\[
  Q_Y\Psi\cong a_*Q_X,
\]
so $\Psi$ is an equivalence. This proves the first isomorphism in
\eqref{eq:reconstruction-kernels}. Taking right adjoints, we obtain
$\Psi G_X\cong G_Y$.
\end{proof}

It remains to determine the action of $\Psi$ on the Mukai lattice.
For $W\in\{X,Y\}$, let
$f_W:=(F_W)_*\colon\widetilde H(\Ku(W),\ZZ)\to K_{\mathrm{top}}(M)$
denote the topological Fourier--Mukai transform, and let
\[
  f_W^H\colon\widetilde H(\Ku(W),\QQ)
  \longrightarrow H^{\mathrm{even}}(M,\QQ)
\]
be the cohomological Fourier--Mukai transform, with the convention of
\cite[\S2.1]{MacriStellari09}.

Let $\iota_W\colon\widetilde H(\Ku(W),\QQ)\hookrightarrow
H^{\mathrm{even}}(W,\QQ)$ be the Mukai-vector embedding
$x\mapsto\ch(x)\sqrt{\td(W)}$, where $\td$ denotes the Todd class.
For the projections $p_W\colon W\times M\to W$ and
$p_M\colon W\times M\to M$, we have
\begin{equation}\label{eq:cohomological-transform}
  f_W^H(x)=(p_M)_*\Bigl(
    p_W^*\iota_W(x)\cup\ch(\cU_W^\vee)
    \cup\sqrt{\td(W\times M)}\Bigr).
\end{equation}
Here we use $F_W\cong\Phi_{\cU_W^\vee[2]}\circ i_W$, as in
Section~\ref{sec:kernel}. Define the universal index
\[
  \Gamma_{\cU_W}(x):=p_{M!}\bigl(p_W^*(x^\vee)[\cU_W]\bigr),
\]
where $p_{M!}$ is the pushforward in topological $K$-theory.
With the dual determinant convention fixed in
Section~\ref{sec:stability}, we have
$\theta_{v_W}(x)=-c_1(\Gamma_{\cU_W}(x))$ on the full topological
orthogonal lattice $v_W^\perp$; see
\cite[Definition~5.4 and Remark~5.5]{BayerMacriProjectivity12} and
\cite[Theorem~29.2 and its proof]{BLMNPSFamilies19}.
By relative two-Calabi--Yau Serre duality, we have
\[
  \bigl(\mathrm Rp_{M*}\RHom(\cU_W,p_W^*A)[2]\bigr)^\vee
  \cong\mathrm Rp_{M*}\RHom(p_W^*A,\cU_W),\qquad A\in\Ku(W).
\]
This is an identity of kernels after projection onto $\Ku(W)$,
so $\Gamma_{\cU_W}(x)=f_W(x)^\vee$ on topological $K$-theory.
Taking first Chern classes and using $c_1(M)=0$, we obtain, after
tensoring with $\QQ$,
\begin{equation}\label{eq:universal-transform-components}
  [f_W^H(x)]_2=c_1(f_W(x))=\theta_{v_W}(x),
  \qquad x\in v_W^\perp\otimes\QQ,
\end{equation}
where $[\cdot]_2$ denotes the degree-two cohomology component.

\begin{proposition}\label{prop:reconstruction-isometry}
The equivalence $\Psi$ in Proposition~\ref{prop:reconstruction}
satisfies $\Psi_*=g$ on the integral Mukai lattices.
\end{proposition}

\begin{proof}
Set $g_\Psi:=\Psi_*$. By $\Psi G_X\cong G_Y$, we have
\[
  \Psi(\cU_{X,m})\cong\cU_{Y,m},
\]
so $g_\Psi(v_X)=v_Y=g(v_X)$.
By the kernel isomorphism $F_Y\Psi\cong F_X$ and functoriality
\cite[\S2.1]{AddingtonThomas14}, we have $f_Y^Hg_\Psi=f_X^H$.
For $x\in v_X^\perp\otimes\QQ$, we obtain, by
\eqref{eq:universal-transform-components} and
\eqref{eq:marking-compatibility},
\[
  \begin{aligned}
    \theta_{v_Y}(g_\Psi(x))
    &=[f_Y^Hg_\Psi(x)]_2=[f_X^H(x)]_2\\
    &=\theta_{v_X}(x)=\theta_{v_Y}(g(x)).
  \end{aligned}
\]
Since $\theta_{v_Y}$ is injective by \eqref{eq:mukai-morphism},
$g_\Psi$ and $g$ agree on $v_X^\perp\otimes\QQ$.
They also agree on $v_X$. Since $v_X^2=2$, we have
\[
  \widetilde H(\Ku(X),\QQ)
  =\QQ v_X\oplus(v_X^\perp\otimes\QQ),
\]
and hence $g_\Psi=g$ over $\QQ$, and over $\ZZ$ by torsion-freeness.
\end{proof}

\begin{corollary}\label{cor:universal-splitting}
There is an isomorphism of Fourier--Mukai kernels
\[
  G_XF_X\cong\id_{\Ku(X)}\oplus[-2]
\]
under which $\eta_X$ is the inclusion of the first summand.
In particular, projection onto the first summand is a retraction of
$\eta_X$.
\end{corollary}

\begin{proof}
Using the kernel identities \eqref{eq:reconstruction-kernels}, we obtain
\[
  G_XF_X\cong\Psi^{-1}\circ(G_YF_Y)\circ\Psi
  \cong\id_{\Ku(X)}\oplus[-2],
\]
where the second isomorphism follows from
Proposition~\ref{prop:universal-splitting}.
These isomorphisms respect the units.
\end{proof}

\begin{remark}
We expect that the analogue of Theorem~\ref{thm:torelli} for
Kuznetsov components of Gushel--Mukai fourfolds can be proved
by a similar argument.
\end{remark}

\section{Derived equivalences of Fano varieties of lines}\label{sec:fano-application}

In this section, we prove Corollary~\ref{cor:fano-ku} using a lattice
cancellation result.
Let $X$ and $Y$ be smooth cubic fourfolds.
For $W\in\{X,Y\}$, set
\[
  \Lambda_W:=\widetilde H(\Ku(W),\ZZ),\qquad
  N_W:=N(\Ku(W)),\qquad T_W:=N_W^\perp\subset\Lambda_W.
\]
Thus $\Lambda_W$ is an even unimodular lattice of signature $(4,20)$,
and $N_W$ is a primitive sublattice of signature $(2,\rk N_W-2)$. For
an even lattice $L$, let $A_L:=L^\vee/L$ denote its discriminant group
and $q_L\colon A_L\to\QQ/2\ZZ$ its discriminant form.

For $W\in\{X,Y\}$, denote by $F(W)$ the Fano variety of lines on $W$.
Let $\Lambda_{\mathrm{Bec}}(F(W))$ be the integral Hodge lattice associated
with $F(W)$ in
\cite[Definition~5.2]{Beckmann23}; it is a derived invariant of $F(W)$
by \cite[Theorem~9.2]{Beckmann23}. By
\cite[Proposition~3.4 and Theorem~0.4]{BottiniHuybrechts25}, there are
Hodge isometries
\begin{equation}\label{eq:fano-lattice-comparison}
  \Lambda_{\mathrm{Bec}}(F(W))\xrightarrow{\ \sim\ }H^2(F(W),\ZZ)\oplus U
  \xrightarrow{\ \sim\ }\Lambda_W\oplus\ZZ\delta_W,\qquad\delta_W^2=-2,
\end{equation}
where $U$ is a hyperbolic plane of type $(1,1)$ and $\delta_W$ is an
algebraic class. These isometries are not canonical, and a Hodge
isometry
$\Lambda_X\oplus\ZZ\delta_X\xrightarrow{\ \sim\ }\Lambda_Y\oplus\ZZ\delta_Y$ need not
send $\delta_X$ to $\pm\delta_Y$. We therefore need a cancellation
result.

\begin{proposition}\label{prop:fano-cancellation}
Let $X$ and $Y$ be smooth cubic fourfolds. Suppose that there is a
Hodge isometry
\[
  \varphi\colon\Lambda_X\oplus\ZZ\delta_X\xrightarrow{\ \sim\ }
  \Lambda_Y\oplus\ZZ\delta_Y,
\]
where $\delta_X,\delta_Y$ are algebraic classes of square $-2$. Then
there is an orientation-preserving Hodge isometry
$\Lambda_X\xrightarrow{\ \sim\ }\Lambda_Y$.
\end{proposition}

\begin{proof}
For $W\in\{X,Y\}$, the discriminant group of
$\Lambda_W\oplus\ZZ\delta_W$ is generated by $\delta_W/2$.
Thus $\varphi(\delta_X)/2=\delta_Y/2+z$ for an algebraic
$z\in\Lambda_Y\oplus\ZZ\delta_Y$.
Consider the orthogonal direct sum
$\Lambda_W\oplus\ZZ\delta_W\oplus\ZZ a_W$,
where $a_W^2=2$ and $a_W$ has Hodge type $(1,1)$. Set
\[
  e_W:=(a_W+\delta_W)/2,\qquad f_W:=(a_W-\delta_W)/2.
\]
These vectors span a hyperbolic plane. The isometry
$\varphi\oplus(a_X\mapsto a_Y)$ maps $e_X$ to $e_Y+z$ and thus
identifies the unimodular overlattices
\begin{equation}\label{eq:fano-stable-U}
  \bigl(\Lambda_W\oplus\ZZ\delta_W\oplus\ZZ a_W\bigr)+\ZZ e_W
  =\Lambda_W\oplus U.
\end{equation}
Restricting this isometry to the algebraic and transcendental parts,
we obtain isometries
\[
  \psi_N\colon N_X\oplus U\xrightarrow{\ \sim\ }N_Y\oplus U,
  \quad \psi_T\colon T_X\xrightarrow{\ \sim\ }T_Y.
\]

We next construct an isometry $\phi_N\colon N_X\xrightarrow{\ \sim\ }N_Y$
inducing the same action on discriminant forms as $\psi_N$.
Set $r:=\rk N_X$, and let $d_1\mid d_2\mid\cdots\mid d_r$ be the
invariant factors of a Gram matrix of $N_X$.
If $r\ge3$, the primitive sublattice $A_2(X)$ has
a basis contained in a basis of $N_X$. The entry $-1$ and the
determinant $3$ of its Gram matrix force $d_1=1$ and $d_2\in\{1,3\}$.
Since $N_X$ is indefinite, it follows from
\cite[Chapter~VIII, Corollary~7.8(1)]{MirandaMorrison09} that it is
determined by its signature and discriminant form and that
$\rO(N_X)\to\rO(q_{N_X})$ is surjective. The same holds for $N_Y$, so we can
choose an isometry inducing the same action on discriminant forms
as $\psi_N$.
When $r=\rk N_X=2$, both lattices are $A_2$, whose discriminant
automorphisms are induced by $\pm\id$.
Let $\gamma_W\colon q_{N_W}\xrightarrow{\sim}-q_{T_W}$ be the
gluing anti-isometry for the unimodular lattice $\Lambda_W$.
By compatibility and gluing, we have
\[
  \bar\psi_T\gamma_X=\gamma_Y\bar\psi_N=\gamma_Y\bar\phi_N,
\]
where $\bar\psi_T$, $\bar\psi_N$, and $\bar\phi_N$ are the induced
isomorphisms of discriminant forms.
Hence $\phi_N\oplus\psi_T$ glues to an integral Hodge isometry
$g_0\colon\Lambda_X\xrightarrow{\ \sim\ }\Lambda_Y$
\cite[\S1.4]{Nikulin80}.
Finally, reflection in $\lambda_{1,Y}$ is an integral Hodge
isometry that fixes the Hodge plane and reverses the orientation of
the positive algebraic plane \cite[Lemma~2.3]{Huybrechts17}.
Composing $g_0$ with this reflection if necessary, we obtain an
orientation-preserving Hodge isometry.
\end{proof}

\begin{proof}[Proof of Corollary~\ref{cor:fano-ku}]
Assume that $\Db(F(X))\simeq\Db(F(Y))$. By
\cite[Theorem~9.2]{Beckmann23}, there is a Hodge isometry
\[
  \Lambda_{\mathrm{Bec}}(F(X))\xrightarrow{\ \sim\ }\Lambda_{\mathrm{Bec}}(F(Y)).
\]
By \eqref{eq:fano-lattice-comparison} and
Proposition~\ref{prop:fano-cancellation}, we obtain an
orientation-preserving Hodge isometry
$\Lambda_X\xrightarrow{\ \sim\ }\Lambda_Y$. The result follows from
Theorem~\ref{thm:torelli}.
\end{proof}

\begin{remark}\label{rem:fano-derived-torelli}
For $W\in\{X,Y\}$, we denote the Markman--Mukai lattice of $F(W)$ by
$\Lambda_{\mathrm{Mar}}(F(W))$ (see \cite[\S9.1.2]{MarkmanSurvey11}).
By \cite[Proposition~7 and Corollary~8]{AddingtonRationality16}, there is a
Hodge isometry between $\Lambda_{\mathrm{Mar}}(F(W))$ and
$\widetilde {H}(\Ku(W),\ZZ)$.
Combining Theorem~\ref{thm:torelli} and
Proposition~\ref{prop:fano-cancellation} with
\cite[Theorem~A]{KemboiSegal25}, we obtain
\[
  \Db(F(X))\simeq\Db(F(Y))
  \quad\Longleftrightarrow\quad
  \begin{gathered}
    \exists\,g\colon\Lambda_{\mathrm{Mar}}(F(X))\xrightarrow{\ \sim\ }
    \Lambda_{\mathrm{Mar}}(F(Y)),\\
    g\text{ an orientation-preserving Hodge isometry}.
  \end{gathered}
\]
This result provides evidence for a derived Torelli principle in
higher-dimensional hyper-K\"ahler geometry.
\end{remark}

\bibliographystyle{amsalpha}
\bibliography{ku_torelli}
\par\addvspace{\baselineskip}

\end{document}